\documentclass[epsf]{article}
\usepackage{amsfonts}
\usepackage{subfigure} 
\usepackage{epsfig}
\usepackage{amsmath, amsthm}
\newcommand{\be}{\begin{eqnarray}}      \newcommand{\ee}{\end{eqnarray}}

\newcommand{\vol}{\mathrm{Vol}}

\newcommand{\rem}{\mathrm{Rm}}
\newcommand{\ric}{\mathrm{Ric}}

\newcommand{\Div}{\mathrm{div}}
\newcommand{\Id}{\mathrm{Id}}
\newcommand{\im}{\mathrm{im}}
\newcommand{\tr}{\mathrm{tr}}
 
\title{Linear and dynamical stability of Ricci flat metrics}
\author{Natasa Sesum} 
\date{} 

\theoremstyle{plain}
\newtheorem{dummy}{Dummy}

\theoremstyle{definition}

\newtheorem{remark}[dummy]{Remark}
\newtheorem{lemma}[dummy]{Lemma}
\newtheorem{theorem}[dummy]{Theorem}
\newtheorem{proposition}[dummy]{Proposition}
  
\newtheorem{definition}[dummy]{Definition}
\newtheorem{step}{Step}[dummy]
\newtheorem{claim}[dummy]{Claim}

\begin{document}

\maketitle

\begin{abstract}
We can talk about two kinds of stability of the Ricci flow at Ricci
flat metrics. One of them is a linear stability, defined with respect
to Perelman's functional $\mathcal{F}$. The other one is a dynamical
stability and it refers to a convergence of a Ricci flow starting at
any metric in a neighbourhood of a considered Ricci flat metric. We
show that dynamical stability implies linear stability. We also show
that a linear stability together with the integrability assumption
imply dynamical stability. As a corollary we get a stability result
for $K3$ surfaces part of which has been done in \cite{dan2002}.
\end{abstract}

\begin{section}{Introduction}

Let $M$ be a closed manifold with a smooth metric $g_0$ on it. We can
flow this metric by the Ricci flow equation
\begin{eqnarray*}
\frac{d}{dt}g_{ij} &=& -2R_{ij}, \\
g_{ij}(0) &=& g_0,
\end{eqnarray*}
that was introduced by R.Hamilton in \cite{hamilton1982}. We can
relate two kinds of stability questions to our initial metric $g_0$.
One of these is just the question of stability of converging Ricci
flows. More precisely, we have a following definition.

\begin{definition} 
Let $g_0$ be a geometry whose Ricci flow $g(t)$ converges. We will say
that $g_0$ is {\it dynamicaly stable} if there exists a neighbourhood
$\mathcal{U}$ of a metric $g_0$ such that the Ricci flow
$\tilde{g}(t)$ of every metric $\tilde{g}\in \mathcal{U}$ exists for
all times $t\in [0,\infty)$ and converges to $g_0$. We will say that
$g_0$ is {\it weakly dynamicaly stable} if there exists a
neighbourhood $\mathcal{U}$ of a metric $g_0$ such that the Ricci flow
$\tilde{g}(t)$ of every metric $\tilde{g}\in \mathcal{U}$ exists for
all times $t\in [0,\infty)$ and converges.
\end{definition}

The other kind of stability is related to Perelman's functional
$\mathcal{F}$ introduced in \cite{perelman2002}. It is given by
$$\mathcal{F}(g,f) = \int_M e^{-f}(|\nabla f|^2 + R)dV_g.$$ 
We can consider functional $\mathcal{\lambda}(g) =
\inf\{\mathcal{F}(g,f)\:\:|\:\:\int_M e^{-f}dV_g = 1\}$. One can compute its 
first variation (\cite{tom2004})
$$\mathcal{D}_g\mathcal{\lambda}(h) = \int_M u\langle -\ric - D^2f,
h\rangle dV_g,$$ 
where $u = e^{-f}$ and we can see that Ricci flat
metrics are critical points of functional $\mathcal{\lambda}$. 
Its second variation at a Ricci flat metric $g$ has been computed in \cite{tom2004} and
$\mathcal{D}^2\mathcal{\lambda}(h,h) = \int_M \langle Lh,h\rangle dV_g$, where
$Lh = \frac{1}{2}\Delta_L h + \Div^*\Div h + \frac{1}{2}d^2v_h$ 
and $\Delta v_h = \Div\Div h $. We can define a linear stability of 
a Ricci flat metric $g_0$ with respect to a second variation of $\mathcal{F}$.
More precisely, 

\begin{definition}
Let $M$ be compact, with $\ric(g_0) = 0$. We will say that $g_0$ is
{\it linearly stable} iff $L \le 0$, that is $\int_M\langle
Lh,h\rangle dV_{g_0} \le 0$ for all directions $h$ (in other words
$L$ has no positive eigenvalues).
\end{definition}

Let $\mathcal{S}_2$ be a bundle of symmetric $2$-tensors on $M$ and
$\mathcal{S}_2^+$ a subset of positive-definite $2$-tensors. Decompose
a tangent space $\mathcal{T}\mathcal{S}_2^+\equiv \mathcal{S}_2$ at a
Ricci flat metric $g_0$ into  subspaces 
\begin{eqnarray*}
\mathcal{S}_2 &=& \im\Div^*\oplus \ker\Div \\
&=& C\oplus E\oplus N\oplus G\oplus S,
\end{eqnarray*}
where
$$C:= \{\delta^*(\delta\eta) : \eta\in \Omega^2\},$$
$$E:= \{\nabla\nabla f:f\in C^{\infty}(M)\},$$
$$N:= \{h\in \mathcal{S}_2 : \delta h = 0, \tr h = 0\},$$
$$S:= \{(\Delta f + \alpha)g - \nabla\nabla f: f\in
C^{\infty}(M,\mathrm{R}),\alpha\in \mathrm{R}\},$$
$$G:= \{\alpha g:\alpha\in \mathrm{R}\}.$$ 
If $M$ is compact and
$\ric(g_0) = 0$, then $L = 0$ on $\im\Div^*$ and on $G$ and $L =
\frac{1}{2}\Delta_L$ on $N$, where $\Delta_L$ is a Lichnerowicz
laplacian $\Delta_l h_{ij} = \Delta h_{ij} + 2R_{ipqj}h^{pq}$,
computed with respect to metric $g_0$.  Since it is already known that
$\Delta_L < 0$ on $C$, $E$, $S$ and $\Delta_L = 0$ on $G$, defining a
linear stability at $g_0$ by saying $L \le 0$ is equivalent to
defining a linear stability at $g_0$ by saying that $\Delta_L \le 0$.

In this paper we want to relate a linear and a dynamical stability at
a fixed Ricci flat metric $g_0$. More precisely, the main theorem we
want to prove in this paper is the following.

\begin{theorem}
\label{theorem-theorem_equivalence}
Let $g_0$ be a Ricci flat metric on a closed manifold $M$. Then if
$g_0$ is dynamicaly stable, it is linearly stable as well. If $g_0$ is
linearly stable and integrable then it is weakly dynamicaly
stable.
\end{theorem}

In section $2$ we will prove the first statement of Theorem
\ref{theorem-theorem_equivalence} by using the monotonicity of
Perelman's functional. In section $3$ we will discuss the
integrability condition of a Ricci flat metric and show how it can be
used to get the convergence of the Ricci flow starting at metrics in a
neighbourhood of a linearly stable and integrable Ricci flat
metric. In section $4$ we will see how Theorem
\ref{theorem-theorem_equivalence} can be applied to get a dynamical
stability of K\"ahler Ricci flat metrics on a $K3$ surface.

{\bf Acknowledgements:} I would like to thank Tom Ilmanen for bringing
the problem of equivalence of those two stabilities to my attention
and for many useful discussions. I would also like to thank ETH in
Z\"urich for its hospitality where a part of this work has been
carried out this summer. I would also like to thank Dan Knopf,
Christine Guenther and Jim Isenberg for discussions about the
stability issue and about their paper \cite{dan2002}.
\end{section}

\begin{section}{From a dynamical stability to a linear stability}

In this section we will prove the following lemma.

\begin{lemma}
\label{lemma-lemma_easy_implication}
Dynamical stability $\Rightarrow$ linear stability.
\end{lemma}

\begin{proof}

We will prove a lemma by contradiction. Assume that $g_0$ is
dynamicaly stable, but that there exists a direction $h_0$ such that
$\int_M\langle Lh,h\rangle dV_{g_0} > 0$. Look at a solution of
\begin{eqnarray*}
\frac{d}{dt}g_{ij} &=& -2R_{ij} \\
g(0) &=& g_0 + h_0.
\end{eqnarray*}
By a dynamical stability assumption we know that there exists
$\lim_{t\to\infty}g(t) = g_0$, since we can choose $h_0$ sufficiently
small. This implies that $\lim_{t\to\infty}\lambda(g(t)) =
\lambda(g_0)$. Let $g_s(t) = g_0 + s(g(t) - g_0)$, $h(t) = g(t) - g_0$
and let $f(s) = \lambda(g_s)$. Then since
$\mathcal{D}_{g_0}\lambda(h(t)) = 0$,
\begin{eqnarray}
\label{equation-equation_ineq}
\lambda(g(t)) - \lambda(g_0) &=& f(1) - f(0) \nonumber \\
&=& \mathcal{D}^2_{g_0}\lambda(h(t),h(t)) + o(|h(t)|^2).
\end{eqnarray}  
Since $\mathcal{D}^2_{g_0}\lambda(h(0),h(0)) > 0$, there exist
$\delta, \eta > 0$ so that $\mathcal{D}^2_{g_0}\lambda(h(t),h(t)) >
\delta |h(t)|^2$ for all $t\in [0,\eta)$. From estimate
(\ref{equation-equation_ineq}) we get that $\lambda(g(t_0)) >
\lambda(g_0) + \epsilon$, for some $\epsilon > 0$ and some $t_0\in
(0,\eta)$. By monotonicity of $\lambda(g(t))$ we have that
$\lambda(g(t)) > \lambda(g_0) + \epsilon$ for all $t\ge t_0$. This is
not possible since $\lim_{t\to\infty}\lambda(g(t)) = \lambda(g_0)$. 
\end{proof} 

\end{section}

\begin{section}{Stability of Ricci flat metrics under the integrability assumption}

In this section we will prove that if a Ricci flat metric is
integrable and linearly stable, then it is dynamicaly stable. The key
step in obtaining this result involves the integrability assumption
which will enable us to get an exponential decay of solutions of a
nonlinear equation whose behaviour is modeled on a behaviour of
solutions of a corresponding linear equation. The Ricci flow PDE
system is not itself strictly parabolic and therefore we will work
with a Ricci DeTurck flow whose PDE system is strictly parabolic and
whose solutions are related to solutions of a Ricci flow by a
$1$-parameter family of diffeomorphisms. We can always find those
gauges for a short time intervals. Our goal is to show that first, we
can extend them all the way up to infinity and second, a solution of a
modified Ricci DeTurck flow has a right kind of behaviour at infinity,
namely an exponential decay. We will do these two things
simultaneously.

We will first give an outline of a proof which involves essentially
five steps. Let $\epsilon'<<\epsilon$ and let $\epsilon$ be very
small.

\begin{enumerate}
\item
Find a gauge (by using DeTurck's trick) in which our solution
$\tilde{g}(t)$ does not differ from $g_0$ more than $\epsilon$ on
some time interval $[0,\eta)$.
\item
Obtain $L^2$ estimates on $|\tilde{g}-g_0|$.
\item
By standard parabolic estimates show that if we fix $A$ arbitrary big,
for sufficiently small initial data we can extend our solution
$\tilde{g}$ to $[0,A)$ such that $|\tilde{g} - g_0|_{C^k} < \epsilon$.
\item
Use the integrability assumption to find new reference Ricci flat
metrics $g_1$ (stationary solutions to the Ricci DeTurck flow) on time
intervals of length $A$, so that we kill zero directions in
$\tilde{g}(t)-g_1$ (zero directions with respect to a Lichnerowicz
laplacian) which yields a decaying type of behaviour.
\item
Use an exponential decay to show that our solution exists for all
times, up to infinity, that we have an exponential convergence of
$\tilde{g}(t)$ and that the same is true for a corresponding solution
of the Ricci flow equation.
\end{enumerate}

Lets first define what the integrability condition means.

\begin{definition}
\label{definition-definition_integrability}
We will say that $g_0$ is integrable, if for any solution {\bf a} of a
linearized deformation equation
$$\mathcal{D}_{g_0}(\ric)({\bf a}) = 0,$$ 
there exists a path $h_u$ of
Ricci flat metrics for $u\in (-\epsilon,\epsilon)$ and $h_0 = g_0$
such that
$$\frac{d}{du}|_{u=0}h_u = a.$$
\end{definition}

In other words, the integrability assumption implies that the set of
metrics $h$ satisfying $\ric(h) = P_{g_0}h = 0$ has a natural smooth
manifold structure near $g_0$.

To finish the proof of Theorem \ref{theorem-theorem_equivalence} we
have to prove the following proposition.

\begin{proposition}
\label{proposition-proposition_dynamic}
Let $g_0$ be a Ricci flat metric on a closed manifold $M$. Assume that
$g_0$ is linearly stable and integrable. Then it is dynamicaly stable.
\end{proposition}

We will give a proof of the Proposition in a sequel of few lemmas and
claims.

\begin{proof}
 Let $h_0$ be a tensor of a small norm, so that $g_0+h_0$ lies
in a small neighbourhood around $g_0$.  More precisely, let
$\epsilon$, $\epsilon'$ be small numbers to be chosen later, so that
$|h_0|_{C^{\infty}} < \epsilon'$, $\epsilon' << \epsilon$.  Fix any
positive integer $k$. Let $\eta$ be such that $\tilde{g}(t)$ exists
and $|\tilde{g}(t)-g_0|_k < \epsilon$ for all $t\in [0,\eta)$. Our
goal is to show that for sufficiently small initial data a solution
$\tilde{g}(t)$ actually exists all the way up to infinity, with
$|\tilde{g}(t)-g_0|_k < \epsilon$ and moreover that it converges
exponentially to $g_0$. Look at the flow
\begin{eqnarray}
\label{equation-equation_flow}
\frac{d}{dt}g_{ij} &=& -2R_{ij} \\
g(0) &=& g_0+h_0 \nonumber. 
\end{eqnarray}
The corresponding modified Ricci DeTurck flow is
\begin{eqnarray}
\label{equation-equation_deturck}
\frac{d}{dt}\tilde{g}(t) &=& -2\ric(\tilde{g}(t)) + P_{g_0}(\tilde{g}), \\
\tilde{g}(0) &=& g_0+h_0 \nonumber.
\end{eqnarray}
Notice that $g_0$ is a stationary solution of Ricci DeTurck flow as
well since $\ric(g_0) = P_{g_0}(g_0) = 0$. 

The linearization of DeTurck operator $\mathcal{A}_{g_0}(g) =
-2\ric(g) + P_{g_0}(g)$ at $g_0$ is
$$\mathcal{D}_{g_0}\mathcal{A}_{g_0}(g))h = \Delta_L h +
\psi_{g_0}h,$$
where $\Delta_L$ is just Lichnerowicz laplacian and
$\psi_{g_0}h = 0$ for a Ricci flat metric $g_0$ (see \cite{dan2002}).
We will use symbols $\Delta_L$ and $L$ interchangebly, to denote a
Lichnerowicz laplacian. We have
\begin{eqnarray}
\label{equation-equation_evolution_h}
\frac{d}{dt}\tilde{h} &=& \frac{d}{dt}\tilde{g} - \frac{d}{dt}g_0 =
\mathcal{A}(\tilde{g}) - \mathcal{A}(g_0) \nonumber \\
&=& \Delta_{g_0}\tilde{h} + U(\tilde{h}) + F(g_0,\tilde{h}),
\end{eqnarray}
and similarly as in \cite{cheeger1994} and \cite{thesis} we have 
a good control over $F$,
\begin{equation}
\label{equation-equation_estimate_F}
|F(\tilde{g},g_0,\tilde{h})|_k \le
C(|\tilde{h}|_k|\nabla^2\tilde{h}|_{k-2} +
|\nabla\tilde{h}|^2_{k-1}),
\end{equation}
where $C$ may depend on bounds on geometries $\tilde{g}(t)$ for $t\in
[0,\eta)$, which is fine since $|\tilde{g}(t) - g_0|_k < \epsilon$ and
therefore geometries of $\tilde{g}(t)$ are uniformly bounded in $t\in
[0,\eta)$ in terms of geometry $g_0$. 

Fix $A>0$ and let $\epsilon > 0$ be small as above. Our goal is to
show that we can extend our solution to $[0,A)$ so that
$|\tilde{h}(t)| < \epsilon$ still holds.  We will first establish some
estimates on $L^2$ norm of $\tilde{h}(t)$ in the following lemma.

\begin{lemma}
\label{lemma-lemma_L2}
There exists a uniform constant $C$ so that for all $t$ for which
$|\tilde{h}|_k < \epsilon$ holds, we have that
$$\int_M|\tilde{h}(t)|^2dV_{g_0} \le e^{C\epsilon A}\int_M|h_0|^2dV_{g_0}.$$
\end{lemma}

\begin{proof}
Remember that $\tilde{h}(t)$ satisfies
$$\frac{d}{dt}\tilde{h}(t) = \Delta_L\tilde{h} + F,$$ 
where $F$ can be controled as in (\ref{equation-equation_estimate_F}).  
Multiply the previous equation by $\tilde{h}$ and integrate it over 
$M$.
\begin{equation}
\label{equation-equation_L2est} 
\frac{d}{dt}\int_M|\tilde{h}(t)|^2dV_{g_0} \le \int_M F\tilde{h}dV_{g_0},
\end{equation}
since $\int_M\langle \Delta_L h,h\rangle dV_{g_0} \le 0$. Notice that we are 
omitting the subscripts. It is not difficult to compute $F$ expilicitely, but
since it is a standard and tedious computation we will omit details here.  
Just to get an idea how the things fit together, we will estimate integrals of
few terms that appear in $F$. Take for example a term
$\nabla\tilde{h}\nabla\tilde{h}\tilde{h}f$, where $f$ is a tensor
depending only on $\tilde{g}$, $g_0$, $\tilde{h}$ and on covariant
derivatives of $g_0$. By partial integration we have
$$\int_M\nabla\tilde{h}\nabla\tilde{h}\tilde{h}f dV_{g_0} =
\int_M(-\tilde{h}\Delta\tilde{h}\tilde{h}f -
\tilde{h}\nabla\tilde{h}\nabla\tilde{h}f  
- \tilde{h}\nabla\tilde{h}\tilde{h}\nabla f)dV_{g_0},$$
and therefore 
$$|\int_M\nabla\tilde{h}\nabla\tilde{h}\tilde{h}f dV_{g_0}| \le
C\int_M|\tilde{h}\Delta\tilde{h}\tilde{h}f dV_{g_0}| + C\int_M
|\tilde{h}|^2|\nabla\tilde{h}\nabla\tilde{f}|dV_{g_0} \le
C\epsilon\int_M|\tilde{h}|^2dV_{g_0},$$
where $C$ basically
depends on bounds on geometry of $g_0$ (notice that $|\tilde{g}|_{C^k}
\le |g_0|_{C^k} + \epsilon \le C + \epsilon$).  Take now a term
$\nabla^2\tilde{h}\tilde{h}f$, where $f$ is of the same form as
before. Then
$$|\int_M\nabla^2\tilde{h}\tilde{h}f\tilde{h}dV_{g_0}| \le
C\epsilon\int_M|\tilde{h}|^2dV_{g_0}.$$
We can estimate other terms that appear in $F$ in a
similar manner and we end up with an estimate 
$$|\int_M F\tilde{h}dV_{g_0}| < C\epsilon\int_M|\tilde{h}|^2dV_{g_0}.$$ 
If we now integrate
(\ref{equation-equation_L2est}) in $t$, since the considered times are 
less or equal than $A$ we get 
$$\int_M|\tilde{h}(t)|^2dV_{g_0} \le e^{CA\epsilon}\int_M|h_0|^2dV_{g_0}.$$
\end{proof}

\begin{lemma}
\label{lemma-lemma_extension}
There exists $\epsilon' = \epsilon'(A,n) << \epsilon$ such that if
$|h_0|_{C^{\infty}} < \epsilon'$, then we can extend $\tilde{g}(t)$
all the way to $[0,A)$ such that $|\tilde{g}(t)-h_0|_k < \epsilon$,
for all $t\in [0,A)$.
\end{lemma}

\begin{proof} 
We will prove that if $\epsilon'$ is sufficiently small, then for
every $t_0\le A$ such that $\tilde{g}(t)$ exists and
$|\tilde{g}(t)-g_0|_k < \epsilon$ on $[0,t_0)$ we actually can improve
our estimates, that is $|\tilde{h}(t)|_k < \epsilon_1$ on $[0,t_1)$
and $\epsilon_1 <<\epsilon$. Therefore, we will be able to extend our
solution past time $t_1$ so that $|\tilde{h}(t)|_k < \epsilon$ holds
past time $t_1$. This will give us that we can extend our solution all
the way up to $A$ so that $|\tilde{h}(t)|_k < \epsilon$ holds on
$[0,A)$.

We will just outline the proof, since it uses standard parabolic
estimates that we carried out in \cite{thesis}. 

\begin{step}
Estimate for $\int_0^{\eta}\int_M|\nabla\tilde{h}|^2dV_{g_0}dt$.
\end{step}

Multiply equation (\ref{equation-equation_evolution_h}) by $\tilde{h}$,
integrate it over $M$ and use the estimates
(\ref{equation-equation_estimate_F}) to get
\begin{equation}
\label{equation-equation_est0}
\frac{1}{2}\frac{d}{dt}\int_M|\tilde{h}|^2 \le
-\int_M|\nabla\tilde{h}|^2 + \int_M\tilde{h}U(\tilde{h}) +
C\int_M(|\tilde{h}|^2|\nabla^2\tilde{h}| +
|\nabla\tilde{h}|^2|\tilde{h}|)dV_{g_0}.
\end{equation} 
Integrate the equation
above in $t\in [0,\eta)$, use the fact that $|\tilde{h}|_k < \epsilon$
for $t\in [0,\eta)$ and use Cauchy-Schwartz inequlity 
$ab \le \theta a^2 + C(\theta)b^2$, where $\theta$ can be taken as small 
as we want to get that
\begin{eqnarray}
\label{equation-equation_est1}
\int_0^{\eta}\int_M|\nabla\tilde{h}|^2dV_{g_0}dt +
\frac{1}{2}\int_M|\tilde{h}(\eta)|^2 &\le&
\frac{1}{2}\int_M|\tilde{h}(0)|^2 +
\theta\int_0^{\eta}\int_M|\nabla\tilde{h}|^2dV_{g_0}dt + \nonumber \\
&+& C\eta\sup_{[0,\eta)}\int_M|\tilde{h}(t)|^2dV_{g_0}.
\end{eqnarray}
By estimate (\ref{equation-equation_est1}) and Lemma
\ref{lemma-lemma_L2} we get that
$\int_0^{\eta}\int_M|\nabla\tilde{h}|^2dV_{g_0}dt$ can be made very
small, comparable to $\epsilon'$ for all $t\in [0,\eta)$.

\begin{step}
\label{step-step_1}
Estimates for
$\int_0^{\eta}\int_M|\frac{d}{dt}\tilde{h}|^2dV_{g_0}dt$,
$\int_0^{\eta}\int_M|\nabla^2\tilde{h}|^2dV_{g_0}dt$ and for
$\sup_{[0,\eta)}\int_M|\nabla\tilde{h}|^2dV_{g_0}$.
\end{step}

\begin{eqnarray}
\label{equation-equation_relation1}
-\int_M\frac{d}{dt}\tilde{h}\Delta\tilde{h} dV_{g_0} &=&
\int_M g_0^{ij}\nabla_i\frac{d}{dt}\tilde{h}\nabla_j\tilde{h}dV_{g_0} \\
&=& \frac{1}{2}\frac{d}{dt}\int_M|\nabla\tilde{h}|^2dV_{g_0}. \nonumber
\end{eqnarray}

\begin{equation}
\label{equation-equation_relation2}
\int_M(\Delta\tilde{h})^2dV_{g_0} = \int_M|\nabla^2\tilde{h}|^2dV_{g_0} + 
\int_M \rem(g_0) * \nabla\tilde{h} * \nabla\tilde{h} dV_{g_0},
\end{equation}
where if $A$ and $B$ are two tensors, we denote by $A*B$ any quantity
obtained from $A\otimes B$ by one of the operations explained in
\cite{dan2004}.

We have the evolution equation for $\tilde{h}$,
\begin{equation}
\label{equation-equation_evolution_h}
\frac{d}{dt}\tilde{h}_{ij} = \Delta\tilde{h}_{ij} +
2R_{ipqj}(g_0)h^{pq} + F(g_0,\tilde{g},\tilde{h}),
\end{equation}
where $\Delta$ is
a laplacian taken with respect to metric $g_0$. Taking the squares of
both sides of equation
$$\frac{d}{dt}\tilde{h}_{ij} - \Delta\tilde{h}_{ij} = 2R_{ipqj}h^{pq}
+ F(g_0,\tilde{g},\tilde{h}),$$
and using relations
(\ref{equation-equation_relation1}) and
(\ref{equation-equation_relation2}) yield 
\begin{equation}
\label{equation-equation_eq1}
\int_M|\frac{d}{dt}\tilde{h}|^2 + \int_M|\nabla^2\tilde{h}|^2 +
\frac{d}{dt}\int_M|\nabla\tilde{h}|^2 \le \int_M(\rem * \tilde{h} +
F(g_0,\tilde{g},\tilde{h}))^2.
\end{equation}
Furthermore,
$$\int_M(\rem * \tilde{h})^2dV_{g_0} \le
C\int_M|\tilde{h}|^2dV_{g_0},$$
\begin{eqnarray}
\label{equation-equation_F}
\int_M|F(g_0,\tilde{g},\tilde{h})|^2dV_{g_0} &\le&
C(\int_M|\nabla^2\tilde{h}|^2|\tilde{h}|^2 +
|\nabla\tilde{h}|^4)dV_{g_0} \\
&\le& C\epsilon\int_M|\nabla^2\tilde{h}|^2dV_{g_0} + 
\epsilon\int_M|\nabla\tilde{h}|^2 dV_{g_0}, \nonumber
\end{eqnarray}
where we can choose $\epsilon$ small enough, so that $C\epsilon <
1/3$.  From Step \ref{step-step_1}, Lemma \ref{lemma-lemma_L2} and
(\ref{equation-equation_F}) it  follows that
$\int_0^{\eta}\int_M|F(g_0,\tilde{g},\tilde{h})|^2dV_{g_0}dt$ can be
made very small, comparable to $\epsilon'<<\epsilon$.  From
(\ref{equation-equation_eq1}) it now follows that
$\int_0^{\eta}\int_M|\frac{d}{dt}\tilde{h}|^2dV_{g_0}dt$,
$\int_0^{\eta}\int_M|\nabla^2\tilde{h}|^2dV_{g_0}dt$ and
$\sup_{[0,\eta)}\int_M|\nabla\tilde{h}|^2dV_{g_0}$ can be made
comparable to $\epsilon'$.

Denote by $H(t) = \frac{d}{dt}\tilde{h}(t)$. If we differentiate
equation (\ref{equation-equation_evolution_h}) in time $t$, we get
\begin{equation}
\label{equation-equation_new}
\frac{d}{dt}H(t) = \Delta_{g_0}H + \rem(g_0)*H +
\frac{d}{dt}(F(g_0,\tilde{g},\tilde{h})).
\end{equation}

\begin{step}
\label{step-step_2}
Estimate on $\sup_{[0,\eta)}\int_M|H|^2dV_{g_0}$ and on
$\int_0^{\eta}\int_M|\nabla H|^2 dV_{g_0}dt$.
\end{step}

If we multiply the equation above by $H$, then integrate it first over
$M$ and then in time $t\in [0,\eta)$, we will get
\begin{eqnarray}
\label{equation-equation_eq2}
\frac{1}{2}\int_MH^2(t)dV_{g_0} + \int_0^{t}\int_M|\nabla H|^2
dV_{g_0} &=& \frac{1}{2}\int_MH^2(0)dV_{g_0} + \int_0^t\int_MH*\rem(g_0)*H
dV_{g_0} + \nonumber \\
&+& \int_0^t\int_M\frac{d}{dt}F HdV_{g_0}.
\end{eqnarray}
We have that
$$\int_M H*\rem(g_0)*H dV_{g_0} \le C\int_M|H|^2dV_{g_0}.$$
It is not
difficult to compute $\frac{d}{dt}(F(g_0,\tilde{g},\tilde{h}))$. We
can estimate each of the terms appearing in it separately, but since
it is a very tedious computation, we will omit details here. Just to
give an idea, some of the terms appearing in
$\frac{d}{dt}(F(g_0,\tilde{g},\tilde{h}))$ are
\begin{eqnarray*}
\int_M\nabla^2H \tilde{h} H f &=& -\int_M\nabla H\nabla\tilde{h}H f
- \int_M \nabla H\tilde{h}\nabla Hf - \int_M \nabla H\tilde{h}H\nabla f \\
&\le& C\epsilon\int_M |\nabla H||H| + C\epsilon\int_M|\nabla H|^2 + 
C\epsilon\int_M|\nabla H||H| \\
&\le& 3\theta\int_M|\nabla H|^2 + C\int_M|H|^2,
\end{eqnarray*}  
where $f$ is a tensor obtained from $\tilde{g}$, $g_0$, $\tilde{h}$
and covariant derivatives of $g_0$ (we can choose $\epsilon$ small, so
that $C\epsilon < \theta$ and $\theta$ is a small positive number).
\begin{eqnarray*}
\int_M\nabla^2\tilde{h} H H f &=& -2\int_M\nabla\tilde{h}\nabla H H f 
-\int_M\nabla\tilde{h}H\nabla f \\
&\le& C\epsilon\int_M |\nabla H|^2 + C\int_M|H|^2 \\
&<& \theta\int_M|\nabla H|^2 + C\int_M|H|^2,
\end{eqnarray*}
where all constants $C$ can be different, but uniform in $t$ and we
will use a same symbol for all of them. We can get similar estimates
for all other terms appearing in
$\frac{d}{dt}(F(g_0,\tilde{g},\tilde{h})$.

All the estimates we have just discussed above, together with
(\ref{equation-equation_eq2}) and results obtained in Step
\ref{step-step_1} and Step \ref{step-step_2} yield that
$\sup_{[0,\eta)}\int_M|H|^2dV_{g_0}$ and $\int_0^{\eta}\int_M|\nabla
H|^2dV_{g_0}dt$ can be made very small, comparable to $\epsilon'$.

If we now consider equation (\ref{equation-equation_new}), using all
the estimates that we have got so far, in the same manner as before we
can get that $\sup_{[0,\eta)}\int_M|\frac{d}{dt}\tilde{h}|^2dV_{g_0}$,
$\int_0^{\eta}\int_M|\nabla^2\frac{d}{dt}\tilde{h}|^2dV_{g_0}dt$ and
$\int_0^{\eta}\int_M(\frac{d^2}{dt^2}\tilde{h})^2dV_{g_0}dt$ can be
made very small, comparable to $\epsilon'$. If we proceed as in the
proof of Proposition $5.1.1$ in \cite{thesis}, by using Sobolev
embedding theorems and standard parabolic regularity theory, we can
get that $|\tilde{h}(t)|_{C^k}$ can be made very small, comparable to
$\epsilon'<<\epsilon$, for all $t\in [0,\eta)$. This means $\tilde{h}$
and therefore $\tilde{g}$ can be extended past time $\eta$ so that
$|\tilde{h}|_k < \epsilon$ still holds past time $\eta$. Actually, we
can start with an arbitrary big $A>0$ so that all the estimates above
depend on $A$ and conclude that we can extend our solution $\tilde{h}$
all the way up to $[0,A)$ so that $|\tilde{h}|_k < \epsilon$ holds.
\end{proof}

We want to show that a solution $\tilde{g}(t)$ exists all the way up
to infinity. If not, let $A' < \infty$ be a maximal real number such
that $\tilde{h}(t)$ exists and $|\tilde{h}(t)|_k < \epsilon$ for all
$t\in [0,A')$. Divide an interval $[0,A')$ by subintervals of length
$A$ and let $N$ be a maximal integer so that $NA < A' <
(N+1)A$. Denote by $I_i = [iA,(i+1)A]$. Over $A_i=M\times I_i$ let
$\pi$ denote orthogonal projection on the subspace
$\ker(-\frac{d}{dt}+\Delta_L)|A_i$ with respect to an inner product
defined by $||\cdot||_{iA,(i+1)A} =
\int_{iA}^{(i+1)A}\int_M|\cdot|dV_{g_0}dt$, where $|\cdot|$ is just a
usual $L^2$ norm. We also have
$$\pi h = (\pi h)_0 + (\pi h)_{\uparrow} + (\pi h)_{\downarrow},$$
where $(\pi h)_0$ represents the radially parallel component
(corresponding to zero eigenvalues of $\Delta L$) and similarly for
$(\pi h)_{\uparrow}$ and $(\pi h)_{\downarrow}$. Since $L \le 0$, $L$
does not have positive eigenvalues and $(\pi h)_{\uparrow} = 0$.  Let
$2\delta = \min\{|\lambda_i| \neq 0\:\:|\:\:\lambda_i$ is an eigevalue
of $\Delta_L\}$. We will briefly describe how we will use the
integrability assumption at this point, with more details following
below. The integrability condition on $g_0$ helps us find a new
stationary solution $g_i$ on every $I_i$ such that
$(\pi(\tilde{g}(t)-g_i))_0 = 0$. This will imply that $\int_M \langle
L(\tilde{g} - g_i), \tilde{g}-g_i\rangle dV_{g_0} \le
-2\delta\int_M|\tilde{g}-g_i|^2dV_{g_0}$ on $I_i$, which yields a
decaying type of behaviour for $\tilde{g}-g_i$.

\begin{lemma}
\label{lemma-lemma_change_reference}
Assume that a Ricci flat metric $g_0$ is integrable in the sense of
Definition \ref{definition-definition_integrability}. If $\alpha =
\alpha(n,A)$ is small enough, then if $\sup_I|\tilde{h}(t)|_k <
\alpha$, where $I$ is one of the intervals $I_i$, there exists a Ricci
flat metric $g_1$ such that $P_{g_0}(g_1) = 0$,
$(\pi(\tilde{g}-g_1))_0 = 0$ and
\begin{equation}
\label{equation-equation_control_g1}
|g_1-g_0| \le C\sup_I |\tilde{g} - g_0|.
\end{equation}
\end{lemma}

\begin{remark}
Before we start proving the Lemma, note that a condition
$P_{g_0}(g_1)=0$ means that a map $\Id:(M,\tilde{g})\to(M,g_0)$ is a
harmonic map.
\end{remark}

\begin{proof}
The proof of Lemma is quite similar to the proof of a corresponding
Lemma $5.56$ in \cite{cheeger1994}, but we will represent it here for
a reader's convenience. The integrability assumption implies that the
set of $g$ satisfying $\ric(g) = P_{g_0}(g) = 0$ has a natural smooth
manifold structure near $g_0$. Let $\mathcal{U}$ be a sufficiently
small euclidean neighbourhood of $g_0$. The tangent space to
$\mathcal{U}$ at $g_0$ is naturally identified with
$$\mathcal{K} = \{a\in \ker\Delta_L \}.
$$ Let $B_i$ be an orthonormal
basis for $\mathcal{K}$ with respect to a natural inner product. Since
$\Delta_L$ is elliptic, a set $\{B_i\}$ is finite. 

Let now $\{\lambda\}$ be a set of eigenvalues for $\Delta_L$ and let
$\{E_{\lambda}\}$ (such that $\Delta E_{\lambda} = -\lambda
E_{\lambda}$) be an orthonormal system of $L^2(\mathcal{S}_2)$ with
respect to a usual inner product. It is then easy to check that
$C_{\lambda} = E_{\lambda}(x)e^{\lambda t}$ is an orthogonal system of
vectors with respect to an inner product $\int_I\langle
\cdot,\cdot\rangle dt$, where $\langle \cdot,\cdot\rangle $ is just a
usual spacelike $L^2$ inner product. If $E_{\lambda}$ is an
eigenvector corresponding to a zero eigenvalue of $\Delta_L$ (there
might be more than one, but finitely many of them, we have denoted
them above by $B_i$), then define $C_{\lambda}' =
\frac{C_{\lambda}}{\sqrt{A}}$; if $E_{\lambda}$ is an eigenvector
corresponding to a nonzero eigenvalue of $\Delta_L$ (at most finitely
many of them) then define $C_{\lambda}'=
\sqrt{\frac{2\lambda}{e^{2\lambda A}(e^{2i\lambda A} -
    1)}}C_{\lambda}$. It is easy to see that $\{C_{\lambda}'\}$ is an
orthonormal system with respect to an inner product defining
$||\cdot||_I := \int_I |\cdot|dt$, where $|\cdot|$ is just a usual
$L^2$ norm. To simplify the notation, write $C_{\lambda}' =
d_{\lambda}e^{\lambda t}E_{\lambda}(x)$, where $d_{\lambda}$ are just
the constants that make $\{C_{\lambda}'\}$ into an orthonormal system.

Define $\psi:\mathcal{U}\to\mathcal{K}$ by $\psi(g) = \sum_i\langle
g,B_i\rangle B_i$. Take any stationary solution of Ricci DeTurck flow
(a metric $g_1$ satisfying $-2\ric(g_1) + P_{g_0}(g_1) = 0$).

\begin{claim}
\label{claim-claim_psi}
$\psi(g_1) = \psi(\pi(g_1)) = (\pi(g_1))_0$.
\end{claim}

\begin{proof}
$g_1 = \sum_{\lambda}a_{\lambda}E_{\lambda}$ and 
\begin{eqnarray*}
\psi(g_1) &=& \sum_i\langle \sum_{\lambda}a_{\lambda} E_{\lambda},B_i\rangle B_i \\
&=& \sum_i a_i B_i,
\end{eqnarray*}
since in the sum $\sum_{\lambda}a_{\lambda} \langle E_{\lambda},B_i\rangle
B_i$ only a term indexed by $i$ survives (remember that
$\{E_{\lambda}\}$ is an orthonormal system with respect to a usual
$L^2$ inner product and $\{B_i\}$ is a subset of that orthonormal
system).  $\{C_{\lambda}'\}$ is a space-time orthonormal basis for a subspace
$\ker(-\frac{d}{dt} + \Delta_L)|_{M\times I}$ and therefore
$$\pi(g_1) = \sum_{\lambda}a_{\lambda}e^{\lambda t}E_{\lambda}(x).$$
Now we have
\begin{eqnarray*}
\psi(\pi(g_1)) &=& \sum_i\langle\sum_{\lambda}a_{\lambda}e^{\lambda t}
E_{\lambda},B_i\rangle B_i \\
&=& \sum_i a_iB_i,
\end{eqnarray*}
since again the only terms that will survive in the sum above will be
those that correspond to a zero eigenvalue of $\Delta_L$ (namely to
$\lambda = 0$, more precisely excatly those for which $E_{\lambda}
= B_i$).  Note that we are using $E_{\lambda}$ to denote all the
orthogonal eigenvectors in our orthonormal system corresponding to the
same eigenvalue $\lambda$. 

We have proved that $\psi(g_1)) = \psi(\pi(g_1))$ and the second term is
$(\pi(g_1))_0$ by a definition of a map $\psi$.
\end{proof} 

Fix some time $t_0\in I$, where $I$ is an interval of length $A$. Note
that a differential of a map $\psi$ is the identity map. By inverse
function therem (keep in mind that our solutions satisfy parabolic
equations and therefore by standard parabolic regularity theory their
$C^k$ norms can be estimated in terms of their $L^2$ norms), when
$|\tilde{g}(t)-g_0|_k$ is small enough on $I$, we get an existence of
$g_1$ such that
$$\psi(g_1) = \sum_i\langle g_1,B_i\rangle B_i =
(\pi\tilde{g}(t_0))_0,$$
which implies $(\pi(g_1-\tilde{g}(t_0)))_0 =
0$, by Claim \ref{claim-claim_psi}. Since $\tilde{g}$ is a solution
and $g_1$ is a stationary solution of the Ricci DeTurck flow we have that
$\tilde{g}(t)-g_1$ satisfies
$$\frac{d}{dt}(\tilde{g}-g_1) = \Delta_L(\tilde{g}-g_1) +
F(g_1,g_0,\tilde{g}-g_1),$$
and it easily follows that
\begin{equation}
\label{equation-equation_all_t}
(\pi(\tilde{g}(t) - g_1))_0 = 0,
\end{equation}
for all $t\in I$. Moreover, from $\psi(g_1) = (\pi\tilde{g})_0$, $\psi(g_0) =
(\pi g_0)_0$ we have that $g_1 = \psi^{-1}((\pi \tilde{g})_0)$ and
$g_0 = \psi^{-1}((\pi g_0))_0$ and therefore
$$||g_1-g_0||_I \le C||\pi^*(\tilde{g}-g_0)||_I.$$
This can be rewritten as
$$|g_1-g_0| \le C\sup_I|\tilde{h}(t)|.$$
\end{proof} 

\begin{lemma}
\label{lemma-lemma_exp_decay}
Let $I=[a,a+A]$ be an interval of length $A$ and let $g_1$ be as in
Lemma \ref{lemma-lemma_change_reference} so that
$(\pi(\tilde{g}(t)-g_1))_0 = 0$ on $I$. Then
$$\sup_{[a+A/2,a+A]}\int_M|\tilde{g}-g_1|^2dV_{g_0} \le
\beta^{-1}\sup_{[a,a+A/2]}\int_M|\tilde{g}-g_1|^2dV_{g_0},$$
where $\beta = e^{\frac{A\delta}{2}}$.
\end{lemma}

\begin{proof}
Notice that $\tilde{g}-g_1$ satisfies
$$\frac{d}{dt}(\tilde{g}-g_1) = \Delta_L(\tilde{g}-g_1) +
F(g_0,g_1,\tilde{g}-g_1),$$
where we can control $F$ as in
(\ref{equation-equation_estimate_F}) by using
(\ref{equation-equation_control_g1}). Denote by $H(t) = \tilde{g}(t)-g_1$.
Similarly as in Lemma \ref{lemma-lemma_L2} we have
\begin{eqnarray*}
\frac{d}{dt}\int_M|H|^2dV_{g_0} &=& \int_M\langle \Delta_L H,H\rangle
dV_{g_0} + \int_M FHdV_{g_0} \\
&\le& -2\delta \int_M|H|^2 + C\epsilon\int_M|H|^2 \\
&\le& -\delta \int_M|H|^2,
\end{eqnarray*}
if we choose $\epsilon$ small so that 
\begin{equation}
\label{equation-equation_choice_eps}
\epsilon < \delta/C.
\end{equation}
This implies that $\int_M|\tilde{g}(t)-g_1|^2dV_{g_0}$ is decreasing
in time $t\in I$ and that
$$\int_M|H(t)|^2dV_{g_0} \le e^{-\delta(t-a)}\int_M|H(a)|^2dV_{g_0}.$$
This yields
\begin{eqnarray*}
\sup_{[a,a+A/2]}\int_M|H(t)|^2 &=& \int_M|H(a)|^2 \\
&\ge& \sup_{[a+A/2,A]}e^{\delta(t-a)}\int_M|H(t)|^2 \\
&\ge& e^{A\delta/2}\sup_{[a+A/2,A]}\int_M|H(t)|^2.
\end{eqnarray*}
\end{proof}

We will now carefully describe what choices for $A$, $\epsilon$ and
$\epsilon'$ we want to make.
\begin{enumerate}
\item
Let $A$ be so big that
\begin{equation}
\label{equation-equation_choice_A}
\frac{1}{1000e^{A\delta}} + \frac{2AC}{e^{A\delta}-1} < e^{-A\delta/4},
\end{equation}
where $C$ is a uniform constant (depending on $g_0$).
\item 
Choose $\epsilon = \min\{\alpha, \delta/C, \delta/(16Cs)\}$,
where $\alpha$ is as in Lemma \ref{lemma-lemma_change_reference} and 
$C$ is some uniform constant that will become apparent later. 
We will discuss $s$ below.
\item
As before, find $\epsilon'$ so small that we can construct 
a solution $\tilde{g}(t)$ to the Ricci DeTurck flow such that 
$\sup_{[0,3A]}|\tilde{h}(t)|_k < \frac{\epsilon}{1000e^{A\delta}}$.
\item 
Note that all the estimates that we have got on
$|\tilde{h}|_{W^{2,k}}$ depend polynomially on $e^{C\epsilon A}$
(this $C$ comes from an estimate in Lemma \ref{lemma-lemma_L2}).
Assume that for the $k$-th order estimate we have a polynomial of degree
$s$.  Therefore, if we choose $\epsilon$ small such that $C\epsilon
< \frac{\delta}{16s}$, then (*) if we start with an initial data
$\tilde{h}(t_0)$, such that $|\tilde{h}(t_0)|_k <
\frac{\epsilon}{e^{A\delta/4}}$ then we will be able to extend
$\tilde{h}(t)$ to an interval $[t_0,t_0+3A)$ so that
$\sup_{[t_0,t_0+3A)}|\tilde{h}(t)|_k < \epsilon$.
\end{enumerate}

Let $A'$ be a maximal real number so that $\tilde{g}(t)$ exists and
$|\tilde{g}(t)-g_0|_k < \epsilon$ for all $t\in [0,A')$. Our goal is
to show that $A'=\infty$. Assume $A' < \infty$ and try to get a
contradiction. Divide an interval $[0,A')$ by subintervals of length
$A$. Call these subintervals $I_i$ as before.

Consider a sequence of sets $A_i = M\times[iA,(i+1)A]$. By Lemma
\ref{lemma-lemma_change_reference}, for every $i$ find $g_i$ such that
for $h_i(t) = \tilde{g}(t)-g_i$ we have $(\pi h_i)_0 = 0$ on
$A_i$. Since $L\le 0$, we have that $\int_M\langle Lh_i,h_i\rangle
dV_{g_0} < 0$ on $I_i$.  Divide an interval $I_i$ by its midpoint,
$A(i+\frac{1}{2})$. By Lemma \ref{lemma-lemma_exp_decay} we have that
$$\sup_{[(i+1/2)A,(i+1)A]}|h_i| \le
\beta^{-1}\sup_{[iA,(i+1/2)A]}|h_i|,$$ 
where $\beta \sim e^{\delta\frac{A}{2}}$ and $\delta$ is such that $2\delta =
\min\{|\lambda_i|\:\:|\:\:\lambda_i\neq 0\}$, where $\lambda_i$ are
the eigenvalues of $\Delta_L$. Since $\tilde{h}(t)$ satisfies 
$$\frac{d}{dt}(\tilde{g}-g_i) = \Delta_L(\tilde{g}-g_i) + F,$$ where
we can control $F$ similarly as in
(\ref{equation-equation_estimate_F}), by Proposition $5.49$ in
\cite{cheeger1994} we get that $\tilde{g}-g_i$ satisfies a decaying
type of behaviour that is
$$\sup_{[(i+1/2)A,(i+1)A]}|h_i| \le\beta^{-1}\sup_{[iA,(i+1/2)A]}|h_i|,$$
implies
$$\sup_{[iA,(i+1/2)A]}|h_i| \le \beta^{-1}\sup_{[(i-1/2)A,iA]}|h_i|.$$
If we apply it inductively we will get that
$$\sup_{[iA,(i+1)A]}|h_i| \le \beta^{-2(i-1)}\sup_{[0,A/2]}|h_i|
= e^{-\delta A(i-1)}\sup_{[0,A/2]}|h_i|.$$
On $M\times I_i$ we have   
\begin{eqnarray*}
|\frac{d}{dt}\tilde{g}| &=& |\frac{d}{dt}(\tilde{g}-g_i)| \\
&\le& |-2\ric(\tilde{g})+2\ric(g_i)| + |P_{g_0}(\tilde{g}) - P_{g_0}(g_i)| \\
&\le& C\sup_{I_i}|\tilde{g}-g_i|_k \le \frac{C\epsilon}{\beta^{2(i-1)}}.
\end{eqnarray*}

\begin{eqnarray*}
\sup_{I_i}|\tilde{g}-g_0| &\le& 2A\sup_{I_i\cup
I_{i-1}}|\frac{d}{dt}\tilde{g}| + \sup_{I_{i-1}}|\tilde{g}-g_0| \\
&\le& 2AC\frac{\epsilon}{\beta^{2(i-2)}} + 2AC\frac{\epsilon}{\beta^{2(i-3)}} + \dots
+ 2AC\frac{\epsilon}{\beta^2} + \sup_{I_3}|\tilde{g}-g_0| \\
&\le& \sup_{I_3}|\tilde{h}| + \frac{2AC\epsilon}{\beta^2 - 1} \\
&\le& \frac{\epsilon}{1000e^{A\delta}} + \frac{2AC\epsilon}{\beta^2 - 1}\\
&\le& \epsilon e^{-A\delta/4},
\end{eqnarray*}
where the last estimate follows by our choice of $A$ in
(\ref{equation-equation_choice_A}. If we take $i=N$, we have that
$\sup_{I_N}|\tilde{g}-g_0|$ can be made smaller than
$\frac{\epsilon}{e^{A\delta/4}}$ and therefore by (*) it can be
extended to $[NA,NA+3A)$ with $\sup_{[NA,(N+3)A)}|\tilde{h}(t)| <
\epsilon$.  Since $A' < (N+3)A $, this contradicts our choice of $A'$.
This means $A' = \infty$ and we can extend our solution $\tilde{g}(t)$
all the way to infinity, so that $|\tilde{g}-g_0| < \epsilon$.
Moreover, we will have
$$|\tilde{g}-g_i| \le Ce^{-\delta t},$$ 
for all $t\in [0,iA)$ and for all $i$. $\{g_i\}$ is a sequence
of Ricci flat metrics with 
$$|g_i - g_0| \le C|\tilde{g} - g_0| \le C\epsilon,$$
and therefore
$\{g_i\}$ have uniformly bounded geometries for all $i$.  This implies
there exists a subsequence of $g_i$ converging to a Ricci flat metric
$g_{\infty}$ satisfying $P_{g_0}(g_{\infty}) = 0$
and
\begin{equation}
\label{equation-equation_exp_conv}
|\tilde{g}(t) - g_{\infty}| < Ce^{-\delta t},
\end{equation}
for all $t\in [0,\infty)$ (notice that we have used a 
standard parabolic regularity
theory to get a $C^k$ convergence, that is we have the following claim).

\begin{claim}
\label{claim-claim_exponential}
We actually have a smooth exponential convergence of $\tilde{g}$ to
$g_0$.
\end{claim}

\begin{proof}
The claim follows from (\ref{equation-equation_exp_conv}) and the
evolution equation (\ref{equation-equation_evolution_h}) for
$\tilde{g} - g_{\infty}$.  We will just sketch a proof here, since
it is pretty standard.  For example, from
(\ref{equation-equation_est0}) it follows that
$\int_M|\nabla\tilde{h}|^2dV_{g_0} \le Ce^{-\delta t}$ for all $t\in
[0,\infty)$. From (\ref{equation-equation_eq1}) it follows that
$\int_M|\nabla^2\tilde{h}|^2dV_{g_0} \le Ce^{-\delta t}$ and
$\int_M|\frac{d}{dt}\tilde{h}|^2dV_{g_0} \le Ce^{-\delta t}$, where
constants $C$ can be different, but we will use a same symbol for
all of them (they are uniform in time $t\in [0,\infty)$). We can
continue a similar analysis and by Sobolev embedding theorems get
that $|\tilde{g}(t) - g_0|_{C^k} \le C(k)e^{-\delta t}$. In other
words, we have an exponential convergence of $\tilde{g}(t)$ towards
$g_0$ in any $C^k$ norm.
\end{proof}

\begin{lemma}
The solution $g(t)$ of the Ricci flow equation converges exponentially
fast to a Ricci flat metric $\bar{g}$.
\end{lemma}

\begin{proof}
$$\frac{d}{dt}g(t) = \ric(g(t)).$$ 
Let $\phi(t)$ be a $1$-parameter
family of diffeomorphisms as in (\ref{equation-equation_exp_conv}). Since the Ricci
tensor is invariant under diffeomorphisms of $M$ we have that 
$\ric(g(t)) = \phi(t)_*\ric(\phi(t)^*g(t))$.
\begin{eqnarray*}
|\ric(\phi(t)^*g(t))|_{k-2} &=& |\ric(\phi(t)^*g(t)-\ric(g_{\infty})|_{k-2} \\
&\le& C|(\phi(t)^*g(t) - g_{\infty}|_k \le Ce^{-t\delta},
\end{eqnarray*}
by (\ref{equation-equation_exp_conv}). This yields 
$|\ric(g(t))|_{k-2} \le Ce^{-t\delta}$ and therefore 
$|\frac{d}{dt}g(t)|_{k-2}\le Ce^{-t\delta}$. This tells us there exists 
some metric $\bar{g}$ so that $|g(t)-\bar{g}|_{k-2}\le Ce^{-t\delta}$. 
$\ric(\bar{g}) = 0$  since
\begin{eqnarray*}
|\ric(\bar{g}| &\le& |\ric(g(t))-\ric(\bar{g})| + |\ric(g(t))| \\
&\le& Ce^{-t\delta} + C|\ric(\phi(t)^*g(t))| \\
&\le& Ce^{-\delta t},
\end{eqnarray*}
and by letting $t\to\infty$ we get that $\bar{g}$ is a Ricci flat 
metric.
\end{proof}
This finishes the proof of Proposition
\ref{proposition-proposition_dynamic}.
\end{proof}

Lemma \ref{lemma-lemma_easy_implication} and Proposition
\ref{proposition-proposition_dynamic} prove Theorem
\ref{theorem-theorem_equivalence}.

If we fix a closed manifold $(M,g)$, we can define another Perelman's
functional
$$\mathcal{W}(g,f,\tau) = (4\pi\tau)^{-n/2}\int_Me^{-f}[\tau(|\nabla
f|^2+R) + f - n]dV,$$
and shrinker entropy by
$$\nu(g) = \inf\{\mathcal{W}(g,f,\tau):f\in
C^{\infty},\:\: \tau > 0,\:\:(4\pi\tau)^{-n/2}\int_M e^{-f}dV = 1\}.$$
In \cite{tom2004} H.D.Cao, R.Hamilton and T.Ilmanen computed the first and the 
second variation of $\nu$. For example, if $(M,g)$ is a positive Einstein
manifold, the second variation $\mathcal{D}^2_g\nu(h,h)$ is given by
$$\frac{d^2}{ds^2}\nu(g(s)) = \frac{\tau}{\vol(g)}\int_M \langle
Nh,h\rangle dV,$$
where 
$$Nh = \frac{1}{2}\Delta_Lh + \Div^*\Div h + \frac{1}{2}D^2v_h -
\frac{g}{2n\tau\vol(g)}\int_M\tr_gh,$$
and
$v_h$ is the unique solution of
$$\Delta v_h + \frac{v_h}{2\tau} = \Div\Div h,$$ 
with $\int_M v_h = 0$. $N$ is degenerate negative elliptic and 
vanishes on $\im\Div^*$. Write 
$$\ker\Div = (\ker\Div)_0\oplus \mathrm{R}g,$$
where $(\ker\Div)_0$ is defined by $\int\tr_g h = 0$. On $(\ker\Div)_0$ 
we have that $N = \frac{1}{2}(\Delta_L - \frac{1}{\tau})$ so an Einstein metric 
$g$ is linearly stable if $\Delta_L \le \frac{1}{\tau}$.

If we now consider a $\tau$-flow $\frac{d}{dt}g = -2\ric(g) +
\frac{1}{\tau}g$ starting at an Einstein metric $g$ with $\ric(g) =
\frac{1}{2\tau}g$, we see that $g$ is a stationary solution of our
flow and therefore it converges to itself when $t\to\infty$. Similarly
as in a Ricci flat case we can get the following theorem.

\begin{theorem}
Let $(M,g)$ be a closed Einstein manifold with Einstein constant
$\frac{1}{2\tau}$, for $\tau > 0$, such that $g$ is linearly stable and
integrable (in the sense as above). Then $g$ is dynamicaly stable.
\end{theorem}

\end{section}

\begin{section}{Stability of K\"ahler Ricci flat metrics on $K3$ surfaces}

\begin{definition}
A $K3$ surface is a connected, closed, smooth, complex surface $M$
with $c_1(M) = 0$ and $b_1(M) = 0$ (that is no global holomorphic
$1$-forms).
\end{definition}

Every $K3$ surface is diffeomorphic to a unique simply-connected
orientable manifold, namely the quartic hypersurface in
$\mathrm{C}\mathrm{P}^3$. Siu proved that every $K3$ surface admits a
K\"ahler metric and by Yau's proof of the Calabi conjecture every
K\"ahler class of a $K3$ surface contains a unique Ricci-flat K\"ahler
metric. We want to fix a K\"ahler Ricci flat metric $g_0$ on $K3$
whose Ricci flow converges to itself and ask whether we get a
convergence of the Ricci flow starting at any metric in a sufficiently
small neighbourhood of $g_0$.

In \cite{dan2002} C.Guenther, J.Isenberg and D.Knopf obtained some
partial results on stability of Ricci flow convergence on a $K3$
surface, using the maximal regularity theory developed by Da Prato and
Gisvard (\cite{prato}) and applied to quasilinear parabolic
reaction-diffusion systems by Simonett (\cite{simonett}) to obtain a
dynamical convergence of the flow to a nontrivial center manifold. Due
to Cao (\cite{cao1985}) every initial K\"ahler metric on a $K3$
surface converges under the Ricci flow to a Ricci flat K\"ahler
metric. We will consider K\"ahler Ricci flat metrics on $K3$.

In \cite{dan2002} it has been shown that for any K\"ahler-Einstein
metric $g_0$ on a $K3$ surface $M$, there is a neighbourhood
$\mathcal{N}_{g_0}$ of $g_0$ in the space of all metrics on $M$ such
that the DeTurck flow $\tilde{g}(t)$ of any initial metric
$\tilde{g}_0$ from $\mathcal{N}_{g_0}$ exponentially approaches a
$58$-dimensional center manifold containing $g_0$, for as long as
$\tilde{g}(t)$ remains in $\mathcal{N}_{g_0}$. C.Guenther, J.Isenberg
and D.Knopf conjectured in \cite{dan2002} that the the Ricci flow of
any initial metric in $\mathcal{N}_{g_0}$ converges to a unique limit
metric in the $58$-dimensional space of K\"ahler Ricci flat metrics
known to exist on a $K3$ surface.

If $(M,g)$ is any Riemannian manifold, denote by $\epsilon(g)$ the
space of infinitesimal deformations of $g$. $h\in \mathcal{S}_2$ is an
infinitesimal Einstein deformation of $g$ if $h\in N$ and
$$\Delta h_{ij} + 2R_{ipqj}h^{pq} = 0.$$ 
If $g$ is a Ricci flat metric
then $\epsilon(g)$ coincides with the kernel of $\Delta_L|N$.

In \cite{dan2002} it was shown that if $(M,g)$ is a K\"ahler-Einstein
metric on a $K3$ surface, then $\Delta_L \le 0$ on $N$ and $\Delta_L <
0$ on $N\backslash\epsilon(g)$. In other words, $\Delta_L$ has no
positive eigenvalues, that is $g$ is linearly stable.

Due to Todorov (\cite{todorov1980},\cite{todorov1983}) it is known
that the infinitesimal deformations of a K\"ahler-Einstein metric
$g_0$ on a $K3$ surface actually correspond to Ricci-flat metrics. In
other words, there is a submanifold
$\mathcal{U}\subset\mathcal{S}_2^+$ of Ricci-flat metrics near $g_0$
such that 
$$\mathrm{T}_{g_0}\mathcal{U} = \epsilon(g_0),$$ 
that is $g_0$ is
integrable in the sense of Definition
\ref{definition-definition_integrability}.

We can now apply Theorem \ref{theorem-theorem_equivalence} to a
linearly stable and integrable K\"ahler Ricci flat metric $g_0$ on a
$K3$ surface $M$ to get the following result.

\begin{theorem}
Let $g_0$ be a K\"ahler Ricci flat metric on a $K3$ surface $M$. There
exists a neighbourhood $\mathcal{N}_{g_0}$ of $g_0$ so that the Ricci
flow of any initial metric in $\mathcal{N}_{g_0}$ converges to a
unique limit Ricci flat metric on $M$, that is $g_0$ is dynamicaly
stable.
\end{theorem}

\end{section}


\begin{thebibliography}{10}

\bibitem{besse} A.Besse: Einstein manifolds; Springer-Verlag,
Heidelberg (1987).

\bibitem{cao1985} H.D.Cao: Deformation of K\"ahler metrics to
K\"ahler-Einstein metrics on compact K\"ahler manifolds;
Invent.Math. 81 (1985) 359--372.

\bibitem{tom2004} H.D.Cao, R.Hamilton, T.Ilmanen: Gaussian densities
for the Ricci flow equation; in preparation.

\bibitem{cheeger1994} J. Cheeger, G. Tian: On the cone structure at
infinity of Ricci flat manifolds with Euclidean volume growth and
quadratic curvature decay; Inventiones Mathematicae 118, (1994),
493--571.
 
\bibitem{dai} Xianzhe Dai, Xiaodong Wang, Guofang Wei: On the
Stability of Riemannian Manifold with Parallel Spinors;
arXiv:math.DG/0311253

\bibitem{dan2004} B.Chow, D.Knopf: The Ricci flow:An Introduction;
American Mathematical Society. 

\bibitem{prato} P.Gisvard, G.Da Parto: Equations d'\'evolution
abstraite nonlin\'eaires de type paraboliques; Ann.Mat.Pura.Appl. 120
(1979) 329--396.

\bibitem{dan2002} C.Guenther, J.Isenberg, D.Knopf: Stability of the
Ricci flow at Ricci-flat metrics; Communications in Analysis and
Geometry.

\bibitem{hamilton1982} R. Hamilton: Three-manifolds with positive
Ricci curvature, Journal of Differential Geometry 17 (1982) 225--306.

\bibitem{perelman2002} G. Perelman: The entropy formula for the Ricci
flow and its geometric applications; arXiv:math.DG/0211159.

\bibitem{thesis} N.Sesum: Limiting behaviour of the Ricci flow; Phd
thesis.

\bibitem{simonett} G.Simonett: Center manifolds for quasilinear
reaction-diffusion systems; Differential Integral Equations 8: 4
(1995) 753--796.

\bibitem{todorov1980} A.Todorov: Applications of the
K\"ahler-Einstein-Calabi-Yau metric to moduli of $K3$ surfaces;
Invent.Math. 61 :3 (1980) 251--265.

\bibitem{todorov1983} A.Todorov: How many K\"ahler metrics has a $K3$
surface?; ''Arithmetic and Geometry:Papers dedicated to
I.R.Shafarevich on the occasion of his sixtieth birthday'', volume 2,
Birkh\"auser (1983) 451--463.

\end{thebibliography}
\end{document}